\documentclass[english,12pt]{amsart}
\usepackage{amssymb}
\usepackage{amsfonts}
\usepackage{amscd}
\usepackage[all]{xypic}

\newcommand{\R}{\mathbb R}

\newtheorem{theorem}{Theorem}
\newtheorem{lemma}{Lemma}
\newtheorem{proposition}{Proposition}

\def\ep{\varepsilon}

\begin{document}

\title{An a priori estimate for the singly periodic solutions of a semilinear equation}
\author{Genevi\`eve Allain}
\address{
Universit\'e Paris-Est Cr\'eteil, Laboratoire d'Analyse et de Math\'ematiques
Appliqu\'ees, UMR CNRS  8050, Facult\'e de Sciences et Technologie,
61, av. du G\'en\'eral de Gaulle, 94010 Cr\'eteil cedex, France} \email{allain@u-pec.fr}

\author{Anne Beaulieu}
\address{
Universit\'e Paris-Est Marne la Vall\'ee, Laboratoire d'Analyse et de Math\'ematiques Appliqu\'ees,
UMR CNRS  8050,
 5 boulevard Descartes,  77454 Marne la Vall\' ee cedex 2, France}\email{anne.beaulieu@univ-mlv.fr}

\date{}
\maketitle

\begin{abstract}
There exists an exponentially decreasing function $f$ such that any singly $2\pi$-periodic positive solution $u$ of $-\Delta u +u-u^p=0$ in 
$[0,2\pi]\times \R^{N-1}$ verifies $u(x_1,x')\leq f(\|Êx'\| )$. We prove that with the same period and with the same function $f$, any singly periodic positive solution of $-\ep^2\Delta u-u+u^p=0$  in $[0,2\pi]\times \R^{N-1}$
verifies $u(x_1,x')\leq f(\|Êx'\| /\ep )$ . We have a similar estimate for the gradient.
\end{abstract}



\section{Introduction.}
Let $N$ be an integer, $N\geq 2$, let $\ep$ and $p$ be positive real numbers, $p>1$. We study the equation
\begin{equation}\label{eq:u}
-\ep^2\Delta u+u-u^p=0\hbox{   in   } S^1\times\R^{N-1}
\end{equation}
where $S^1=[0,2\pi]$. We mean that $u$ is $2\pi$-periodic in $x_1$.
We consider the positive solutions of (\ref{eq:u}), $u(x_1,x')$ ($x_1\in S^1$ and $x'\in\R^{N-1}$) that tend to 0 as $\| x'\|$ tends to $\infty$, uniformly in $x_1$. It is known that these solutions are radial in $x'$ and decreasing in $\| x'\|$.  This can be proved
by an application of the moving plane method (\cite{Gidas},
\cite{Almeida},  \cite{Berestycki}). The ground-state solution $w_0$, defined and radial on $\R^{N-1}$ is a particular solution which does not depend on $x_1$.
In \cite{Dancer}, Dancer proved the existence of positive solutions really depending on $x_1$ and $x'$. In \cite{singly}, we studied the case $N=2$ and we proved the following result:
 \begin{theorem} \label{theoremsingly}  (i) The first continuum $\Sigma_1$ of positive bounded solutions even in $x_1$ and $x'$ of
(\ref{eq:u}) bifurcating from
$(\ep_{\star},w_0({x'/\ep_{\star}}))$ is composed of $(\ep_{\star},w_0({x'/\ep_{\star}}))$ and of all the
solutions $(\ep,z)$ of $(\ref{eq:u})$ such that $z>0$, $z$ even in $x_1$ and
$x_2$, $\lim_{x_2\rightarrow\infty}z=0$ and ${\partial
z\over\partial x_1}<0$ in $]0,\pi[\times\R+$. \\
(ii)There exists a bounded
subset $\mathcal A$ of $L^{\infty}(S^1\times\R+)$ such that the set
$\Sigma_1$ is entirely contained in $]0,\ep_{\star}]\times{\mathcal A}
$.\\
(iii) For each $(\ep,z)\in \Sigma_1$, $z$ is an isolated point of 
$\{ v\in L^\infty (S^1\times \R+ );$ $v$  even in $x_1$ and $x'$;
 $(\ep, v)$ solution of $(\ref{eq:u}) \}$.
For every $\ep>0$, $\ep<\ep_{\star}$, there exists a finite
number of solutions $(\ep,z)$ in $\Sigma_1$. \\
(iv) There exists $\ep_0$
such that for all $0<\ep<\ep_0$ this continuum is a curve that has a
one to one ${\mathcal C}^1$ parameterization $\ep\rightarrow
(\ep,z_{\ep})$. 
\end{theorem}
In this paper we suppose that
\begin{equation}\label{eq:souscritique}
1<p<{N+2\over N-2}
\end{equation}
If $N=2$, this condition is $p>1$. 

We know, see \cite{Dancer}, that the condition (\ref{eq:souscritique}) for $p$ is a necessary and sufficient condition to have the following property :
There exists $M>0$ such for all $\ep>0$, any positive solution $u$ of (\ref{eq:u}) verifies 
\begin{equation}\label{bornŽ}
\|u\|_{L^{\infty}}\leq M
\end{equation}
This property is related to the nonexistence of positive  solutions for the equation $-\Delta u-u^p=0$, more precisely
\begin{equation}\label{eq:GS}
( v\geq0\quad,\quad -\Delta v-v^p=0\hbox{  in  }\R^N)\Rightarrow (v=0)
\end{equation}
(see Gidas and Spruck \cite{GidasSpruck}).
This paper is devoluted to some a priori  estimates for the solutions of (\ref{eq:u}).

\begin{theorem} \label{exponentiel} There exists a real number $K$ independent of  $\ep>0$ and of any solution $u$ of (\ref{eq:u}), such that for all $x=(x_1,x')$ in $S^1\times \R^{N-1}$, we have, with  $r'=\|x'\|$, 
\begin{equation} \label{eq:uasympt}
u(x)\leq K e^{-\frac{r'}{\ep}}\Big(\frac{r'}{\ep} \Big)^{2-N\over 2} 
\end{equation}
\begin{equation} \label{eq:gradasympt}
\|\nabla u(x)\|\leq  {\frac{K}{\ep}}e^{-\frac{r'}{\ep}}\Big(\frac{r'}{\ep} \Big)^{2-N\over2}
\end{equation}
\end{theorem}

In \cite{singly}, we have proved (\ref{eq:uasympt}) for $N=2$ but with a constant  $K$ depending on the solution $(\ep, u)$ for $\ep$ greater than some $\overline \ep >0$. Our proof extends easily for $N\geq 2$ and for the derivatives of $u$. We have now to prove that $K$ is independent from the solution $(\ep,u)$, even when $\ep$ tends to 0.\\

In all what follows we will use $\tilde u(x_1,x')=u(\ep x_1, \ep x')$ for $(x_1,x')\in S^1/\ep \times \R^{N-1}$. The notation $\Delta'$ will stand for the Laplacian operator in $\R^{N-1}$.

\section{Proof of Theorem \ref{exponentiel}} We begin the proof by two propositions. 
\begin{proposition}\label{GG} Let $v$ be a bounded solution of 
\begin{equation}\label{eq:v}
-\Delta v+v-v^p=0\hbox{ in  }\R^2
\end{equation} 
 Let us suppose that ${\partial v\over\partial x_i}$ is bounded and  ${\partial v\over\partial x_i}\leq0$ in $\R^2$, for  $i=1$ or $i=2$. Then $v$ does not depend on the variable $x_i$.
\end{proposition}
{\bf Proof: } 
From  Ghoussoub-Gui, \cite{GG}, Theorem 1.1, there exist a function $U$ and a vector $a\in\R^2$ such that $v(x)=U(a.x)$. We have 
$$-\|a\|^2U''(a.x)+(U-U^p)(a.x)=0$$
If $a_i\neq 0$,  $U$ is monotone and bounded  in $\R$. The only possibility is that $U$ is a constant function, equal to 0 or 1, so is $v$.

\begin{proposition}\label{limuniformegeneral} Let $(\ep,u)$ be solutions of (\ref{eq:u}). Then $  \tilde u(x_1,r)$
tends to 0 as $r$ tends to $\infty$, uniformly with respect to
$x_1$ and to $\ep$ and $u$.
\end{proposition}

{\bf Proof} 
In this proof, we omit the indices of the sequences.  Let us suppose, by contradiction, that there exist a sequence $(a,b)\in \R^N$,
with $\|b\|$ tending to $+\infty$,  a real positive number $\ep_2$ and solutions $(\ep,u)$ of (\ref{eq:u})
such that $\tilde  u(a,b)\geq \ep_2$. We can suppose $\ep_2<1$. For every solution $(\ep,u)$, we have
that $\lim_{r\rightarrow\infty}\tilde  u(x_1,r)=0$, uniformly in $x_1$. So, for
every $\ep_1\in]0,\ep_2[$, there exists a sequence, $\overline
b$, $\|\overline b\|\geq \| b\|$, such that  $\tilde  u(a,\overline b)=
\ep_1$.  With the same argument, we define a
sequence, still denoted by $b$, with $\|b\|$ tending to $\infty$, such that
$\tilde  u(a,b)= \ep_2$. As $\tilde u$ is radial in $x'$, let us define
$$v(x_1,r)=\tilde  u(x_1+a,r+\| b\|) \hbox{ for } r\geq - \| b\|$$
and
$$\overline v(x_1,r)=\tilde  u(x_1+a, r+\|\overline b\|) \hbox{ for } r\geq - \| \overline b\|$$
The function $v$ verifies
$$ -v_{x_1x_1}-v_{rr}-\frac{N-2}{r+\| b\|}v_r+v-v^p=0$$
and $\overline v$ verifies a similar equation.
It is standard that the both sequences $v$ and $\overline v$
tend uniformly on the compact sets of $\R^2$ to limits, which will
be denoted respectively by $z$ and $\overline z$. But $z$ and
$\overline z$ are positive, bounded and non increasing in the variable $r$ and they are periodic in $x_1$. Moreover, $z$ and $\overline z$ verify
$$ -z_{x_1x_1}-z_{rr}+z-z^p=0 \hbox{ in } \R^2$$
 By Proposition \ref{GG}, $z$ and $\overline z$ depend
only on $x_1$.
 By Kwong, \cite{Kwong}, if they are not constant functions, they oscillate
indefinitely as $x_1$ tends to $\infty$, around the solution $1$. As $0<\ep_1<\ep_2<1$, then $\overline z$ and $z$ are not constant solutions.   So $z$ and $\overline z$ oscillate
infinitely around $1$, too. 
The function $h=z-1$ and the function $\overline h=\overline z-1$ verify respectively the equations
\begin{equation}\label{eq:z}
h''+h(-1+{z^p-1\over z-1})=0\quad\hbox{  and  }\quad\overline h''+\overline h(-1+{ \overline z^p-1\over \overline z-1})=0
\end{equation}
As $z\geq \overline z$ and $z(0)>\overline z(0)$, we have from the ordinary differential equations theory that $z>\overline z$. It is easy to see that $-1+{z^p-1\over z-1}>-1+{\overline z^p-1\over
\overline z-1}$. By the Sturm Theory (see Ince, quoted in \cite{Kwong} , Lemma 1), applied to the equations
(\ref{eq:z}), there exists at least a zero of
$z-1$ between any two consecutive zeroes of $\overline z-1$. But
there exist pairs $(\alpha,\beta)$ of zeroes of $ \overline z-1$ such
that $\overline z>1$ in $]\alpha,\beta[$. Thus $z>1$ in
$[\alpha,\beta]$. We get a contradiction. 
We infer that
the sequence $(a,b)$, in the beginning of this proof, doesn't
exist. We have proved the proposition.\\

We will need the following lemma

\begin{lemma}\label{gradientborne}
There exists $M$, such that for all solution $(\ep,u)$ of (\ref{eq:u})
\begin{equation} \label{eq:bddgrad}
\|\nabla \tilde u\|_{L^{\infty}({S^1\over\ep}\times\R^{N-1})}\leq M
\end{equation}
\end{lemma} 
{\bf Proof}  Let  $(a,b)\in(S^1/\ep)\times \R^{N-1}$. We set $v(x_1,x')=\tilde u(x_1+a,x'+b)$. It verifies  $-\Delta v+v-v^p=0$ in $(S^1/\ep)\times\R^{N-1}$. Moreover, we have $\|v\|_{\infty}\leq M$, for a constant $M$ independent from $\ep$. By standart elliptic arguments, \cite{Gilbarg}, $\nabla v$ is bounded on the compact sets of $\R^N$. So, there exists $M$, independent from $\ep$, such that $\|\nabla v(0,0)\|\leq M$. This proves (\ref{eq:bddgrad}).\\

{\bf Proof of Theorem \ref{exponentiel}}  We define
$$h(r')=\int_0^{2\pi/\ep}\tilde u(x_1,r')dx_1$$
There exists a constant $C$, independent from the solution $(\ep,u)$, such that
$\|h\|_{L^{\infty}(\R^{N-1})}\leq {C\over\ep}$.
Since $u\rightarrow0$, uniformly in $x_1$, as $r'$ tends to $\infty$, then
for all $\eta<1$, there exists $X>0$
such that for all $ r'>X$ and for all $\ep $ we have for all solution $(\ep,u)$ and for all $x_1\in{S^1\over\ep}$
\begin{equation}\label{eq:up-1}
\tilde u^{p-1}(x_1,r')<\eta
\end{equation}
Integrating (\ref{eq:u}) with respect to $x_1$, we find for $r'>X$
\begin{equation}\label{eq:hr} 
h_{rr}+((N-2)/r')h_r>(1-\eta)h
\end{equation}
Let us multiply (\ref{eq:hr}) by $h_r$, we obtain that the function $h_r^2-(1-\eta)h^2$ is non increasing. Moreover, it tends to 0 as $r'$ tends to $\infty$. We get $h_r+\sqrt{1-\eta}h\leq0$, for $r'>X$. So there exists $C$ such that for all $r'$
\begin{equation}\label{h}
h(r')\leq {C\over\ep}e^{-\sqrt{ 1-\eta}r'} 
\end{equation}
Let us remark that the constant $C$ is independent from the choice of the solution $(\ep,u)$.\\

Let $R>0$ be a given positive
real number. We use a Harnack inequality (\cite{Gilbarg}, Theorem 9.20) to get that
there exists a constant $C$ independent from $y$ and
from $\ep$ such that
\begin{equation}\label{eq:harnack}
\sup_{B_R(y)}{\tilde u}\leq
C\int_{B_{2R}(y)}{\tilde u}\leq C\int_{\|x'-y'\|\leq 2R}\int_{y_1-R}^{y_1+R}\tilde u(x_1,x')dx_1dx'
\end{equation}
that gives
$$\sup_{B_R(y)}{ \tilde u}\leq C\int_{\|x'-y'\|\leq 2R}h(\|x'\|)dx'.$$
Finally, using (\ref{h}), for all $\eta\in]0,1[$ there exists $C$, independent from the solution $(\ep,u)$, such that,
\begin{equation}\label{eq:expeta}
\tilde u(y)\leq {C\over\ep}e^{-\eta\|y'\|}
\end{equation}
For the remainder of the proof, we will need the Green function for the equation (\ref{eq:u}). We have
\begin{equation}\label{eq:fonction de Green}
G(x_1,x')=\sum_{j=0}^{\infty}{k_j^{N-3}\over\ep^{N-1}}g({k_j\over\ep}x')\cos(jx_1)
\end{equation}
where $k_j=\sqrt{1+\ep^2j^2}$ 
and $g$ is the Green function for the operator $-\Delta' +I$ in $\R^{n}$, $n=N-1$, with the null limit at infinity. It is recalled in \cite{Gidas} that
\begin{equation}\label{eq:bessel}
0<g(r)\leq C{e^{-r}\over r^{n-2}}(1+r)^{(n-3)/2}\hbox{  for  } n\geq 2\hbox{  and  } g(r)={1\over2}e^{-r}\hbox{ for }n=1
\end{equation}
We will need the following estimate, valid for all $\eta\in]0,1[$.
\begin{equation}\label{eq:integraleeta}
\int_{\R^{N-1}}g(\|y'-x'\|)e^{-\eta\|y'\|}dy'\leq Ce^{-\eta\|x'\|}
\end{equation}
which is an easy consequence of (\ref{eq:bessel}).
For all function $f$, that is $2\pi$-periodic in $x_1$, the solution of
$$-\ep^2\Delta u+u=f\hbox{  in  }\R^N$$
that is $2\pi$-periodic in $x_1$ and that tends to 0, as $\|x'\|$ tends to $\infty$ is
$u=G\star f$.
If $f$ is positive, then $u$ is positive, by the maximum principle. So $G$ is positive.
Moreover we can use (\ref{eq:fonction de Green}) to verify that
\begin{equation}\label{eq:integraleG}
\int_{S^1}G(x_1,x')dx_1={2\pi\over\ep^{N-1}} g({x'\over\ep})
\end{equation}

Let us prove that for all $\eta\in]0,1[$, there exists $C$, independent from $x_1$ and from $(\ep,u)$ such that 
\begin{equation}\label{eq:ueta}
\tilde u(x_1,x')\leq {C}e^{-\eta r'}
\end{equation}
It is clear by (\ref{eq:expeta}) that for all solution $(\ep,u)$ and all $\eta\in]0,1[$, the function 
$\tilde u e^{\eta r}$ belongs to $L^{\infty}(\R^N)$.  We set 
$$K(\eta)=\|\tilde u e^{\eta r'} \|_{\infty}$$
We use the Green function $G$ to get
\begin{equation}\label{eq:uetup}
u(x_1,x')=\int_{S^1\times\R^{N-1}}G(y_1-x_1,y'-x')u^p(y_1,y')dy_1dy'
\end{equation}
and (\ref{eq:integraleG}) gives
$$\tilde u(x_1,x')\leq 2\pi K\Big(\frac{\eta}{ p}\Big)^p\int_{\R^{N-1}}g(\|y'-x'\|)e^{-\eta\|y'\|}dy'$$
By (\ref{eq:integraleeta}), we infer that there exists a constant $C$, independent from $(\ep,u)$, such that
\begin{equation}\label{eq:Ketasurp}
K(\eta)\leq C K\Big(\frac{\eta}{ p}\Big)^p
\end{equation}
Now, let $\tau=(\tau_1,\tau')$ be such that the function $\tilde u( x+\tau)e^{\eta\|x'+\tau'\|}$ attains its maximal value at $x=0$. The existence of $\tau$ is provided by (\ref{eq:expeta}). 
Let us suppose that $K(\eta)$ tends to $\infty$. We claim that $\|\tau'\|$ tends to infinity. Let us prove this claim. Let $\alpha$ be a positive real number, that will be chosen later.  We set
$$v(x)=\tilde u(\alpha x+\tau)e^{\eta\|\alpha x'+\tau'\|}/K(\eta)$$
It verifies
$$-\Delta v+\Big(1+\eta^2+\frac{(N-2)\eta}{\|\alpha x'+\tau'\|}\Big)\alpha^2v$$

$$=K(\eta )^{p-1}\alpha^2e^{(-p+1)\eta\|\alpha x'+\tau'\|}v^p+{2\eta\alpha^2\over K(\eta)}\sum_{i=2}^N {\partial \tilde u \over \partial x_i}(\alpha x+\tau ){(\alpha x_i+\tau_i)\over \| \alpha x'+\tau'\|}
e^{\eta\|\alpha x'+\tau'\|}$$
If $\|\tau'\|$ were bounded, we would choose $\alpha$ that tends to 0 such that $K(\eta)^{p-1}\alpha^2e^{(-p+1)\eta\|\alpha x'+\tau'\|}$ tends to 1.  By Lemma \ref{gradientborne} and
by standard results, $v$ would tend to a limit $\overline v$, uniformly in the compact sets of $\R^N$. Then, $\overline v$ would verify
$-\Delta\overline v-\overline v^p=0$
while $0\leq\overline v\leq 1$ and $\overline v(0)=1$. This is impossible by (\ref{eq:GS}). So, if we suppose that $K(\eta)$ tends to $\infty$, then $\|\tau'\|$ tends to $\infty$. Let $\tilde \tau=(\tilde \tau_1,\tilde\tau')$ be such that $K(\frac{\eta}{ p})=\tilde u(\tilde\tau)e^{{\eta\over p}\|\tilde\tau'\|}$. We have
$K(\frac{\eta}{ p})^p=\tilde u^p(\tilde\tau)e^{\eta\|\tilde\tau'\|}$, that gives
\begin{equation}\label{eq:tildetau}
K\Big(\frac{\eta}{ p}\Big)^p\leq K(\eta)\tilde u^{p-1}(\tilde\tau)
\end{equation}
Then (\ref{eq:Ketasurp}) and (\ref{eq:tildetau}) give
\begin{equation}\label{eq:doubleinegalite}
K(\eta)\leq C K\Big(\frac{\eta}{ p}\Big)^p\leq C K(\eta)\tilde u^{p-1}(\tilde\tau)
\end{equation}
Consequently, if $K(\eta)$ tends to $\infty$, then $K(\frac{\eta}{ p})$ tends to $\infty$, too. Then, $\|\tilde\tau'\|\rightarrow\infty$. By Proposition \ref{limuniformegeneral}, we have $\tilde u(\tilde\tau)\rightarrow0$. Then (\ref{eq:doubleinegalite}) gives a contradiction. So, we have proved that for all $\eta\in]0,1[$, $K(\eta)$ is bounded, independently from $(\ep,u)$. We have (\ref{eq:ueta}). Now, let us choose $\eta$ such that $\eta p>1$. In \cite{Gidas}, it is proved that for $b>1$  and for $N-1\geq2$
\begin{equation}\label{eq:b1}
\int_{\R^{N-1}}g(\|x'-y'\|)e^{-b\|y'\|}dy'\leq C \|x'\|^{2-N\over2}e^{-\|x'\|}
\end{equation}
We can use (\ref{eq:bessel}) to prove that the estimate (\ref{eq:b1}) is valid also for $N=2$. Now we use (\ref{eq:uetup}), (\ref{eq:integraleG}) and (\ref{eq:b1}) to obtain (\ref{eq:uasympt}) with $K$ independent from $(\ep, u)$.
Now, let us estimate the gradient of $u$.
We have, for $i=1,...,N$
\begin{equation}\label{eq:dusurdx}
{\partial u\over \partial x_i}(x_1,x')=p\int_{S^1\times\R^{N-1}}G(y_1-x_1,y'-x')(u^{p-1}{\partial u\over \partial x_i})(y_1,y')dy_1dy'
\end{equation}
Since ${\partial\tilde u\over \partial x_i}$ is bounded and $u\leq C e^{- r'/\ep}$, that gives 
$$\mid {\partial u\over \partial x_i}(x_1,x')\mid\leq {C\over \ep}\int_{S^1\times\R^{N-1}}G(y_1-x_1,y'-x')e^{-(p-1)\|y'\|/\ep}dy_1dy'$$
and (\ref{eq:integraleG}) gives
\begin{equation}\label{eq:partial}
\mid {\partial \tilde u\over \partial x_i}(x_1,x')\mid\leq C\int_{\R^{N-1}}g(\|y'-x'\|)e^{-(p-1)\|y'\|}dy'
\end{equation}
Now, the proof is more easy if $p>2$ than if $p<2$. If $p>2$, we deduce directly (\ref{eq:gradasympt}) from (\ref{eq:b1}) and (\ref{eq:partial}). If $1<p<2$, we deduce from (\ref{eq:integraleeta}) and (\ref{eq:partial}) that
$$\mid {\partial \tilde u\over \partial x_i}(x_1,x')\mid\leq Ce^{-(p-1)\|x'\|}$$
Iterating this process, we get an integer $k$ such that
$$\mid {\partial \tilde u\over \partial x_i}(x_1,x')\mid\leq Ce^{-k(p-1)\|x'\|}$$
with $k(p-1)<1$ and $(k+1)(p-1)\geq1$. If $(k+1)(p-1)>1$, we get (\ref{eq:gradasympt}) and the proof is complete. If $(k+1)(p-1)=1$, we get
\begin{equation}\label{eq:partial2}
\mid {\partial\tilde u\over \partial x_i}(x_1,x')\mid\leq C\int_{\R^{N-1}}g(\|y'-x'\|)e^{-\|y'\|}dy'
\end{equation}
If $N\geq 3$, we have
$$\int_{\R^{N-1}}g(\|y'-x'\|)e^{-\|y'\|}dy'\leq C\int_{\R^{N-1}}e^{-\|x'-y'\|-\|y'\|}(1+\|x'-y'\|)^{(N-4)/2}/\|x'-y'\|^{N-3}dy'$$
We can write the integral in the right hand member of this inequality as $I=I_1+I_2$ and
$$I_1=\int_{\|z\|\leq\|x'\|}e^{-\|z\|-\|z+x'\|}(1+\|z\|)^{(N-4)/2}/\|z\|^{N-3}dz$$
and 
$$I_2=\int_{\|z\|\geq\|x'\|}e^{-\|z\|-\|z+x'\|}(1+\|z\|)^{(N-4)/2}/\|z\|^{N-3}dz$$
We obtain, as $\|x'\|$ tends to $\infty$,
$$I_1\leq e^{-\|x'\|}\int_0^{\|x'\|} (1+s)^{N-4\over2}sds\quad
\hbox{ and  }\quad
I_2\leq e^{\|x'\|}\int_{\|x'\|}^{+\infty}e^{-2s}(1+s)^{N-4\over2}sds$$
These integrals are both less than $C e^{-\|x'\|}\|x'\|^{N\over2}$.
Thus, if $N\geq3$ we have obtained that
\begin{equation}\label{eq:p2}
\mid {\partial\tilde u\over \partial x_i}(x_1,x')\mid\leq C e^{-\|x'\|}\|x'\|^{N\over2}
\end{equation}
 If $N=2$, we have, when $|x'|$ tends to $\infty$
$$\int_{\R}g(|y'-x'|)e^{-|y'|}dy'\leq C\int_{\R}e^{-|x'-y'|-|y'|}dy'\leq C|x'|e^{-|x'|}$$
In any case, we get that there exists $b\in]0,1[$, with $b+p-1>1$ and such that
$\mid {\partial\tilde u\over \partial x_i}(x_1,x')\mid\leq C e^{-b\|x'\|}$
Using this estimate in (\ref{eq:dusurdx}) and thanks to (\ref{eq:b1}), we get (\ref{eq:gradasympt}), for $1<p<2$. If $p=2$, (\ref{eq:partial}) is (\ref{eq:partial2}) and we deduce (\ref{eq:p2}) again. This ended the proof of Proposition 
\ref{exponentiel}.

\end{document}